\documentclass[12pt]{article}
\usepackage{latexsym,amssymb,amsmath,mathrsfs}
\usepackage[anchorcolor=blue,%
 bookmarks=true,%
 bookmarksnumbered=true,%
 dvipdfm,%
]{hyperref}
\makeatletter \@addtoreset{equation}{section} \makeatother

\newcommand{\nn}{\nonumber}

\newenvironment{demo}[1]{\par\begin{trivlist}%
\item[]{\bf #1}\ }{\end{trivlist}\par}
\newcommand{\Proof}{\begin{demo}{\bf Proof.\ }}
\newcommand{\Proofof}[1]{\begin{demo}{\bf Proof of #1.\ }}
\newcommand{\qed}{\hfill $\langle {\bf q.e.d.}\rangle$ \end{demo}}
\newcommand{\toy}{\ \rule[0em]{0.5ex}{1.8ex}}
\newcommand{\QED}{\toy\end{demo}}
\newcommand{\eop}{\hfill $\square$ \end{demo}}

\newtheorem{thm}{Theorem}[section]
\newtheorem{prop}[thm]{Proposition}
\newtheorem{cor}[thm]{Corollary}
\newtheorem{lem}[thm]{Lemma}
\newtheorem{defn}[thm]{Definition}
\newtheorem{remark}[thm]{Remark}
\newtheorem{example}[thm]{Example}
\newtheorem{assumption}{Assumption}

\newcommand{\bthm}[1]{\begin{thm}\label{th:#1}}
\newcommand{\bthmn}[2]{\begin{thm}{\bf (#2)}\label{th:#1}}
\newcommand{\bprop}[1]{\begin{prop}\label{prop:#1}}
\newcommand{\bpropn}[2]{\begin{prop}{\bf (#2)}\label{prop:#1}}
\newcommand{\blem}[1]{\begin{lem}\label{lem:#1}}
\newcommand{\blemn}[2]{\begin{lem}{\bf (#2)}\label{lem:#1}}
\newcommand{\bcor}[1]{\begin{cor}\label{cor:#1}}
\newcommand{\bdefn}[1]{\begin{defn}\label{def:#1}\rm}
\newcommand{\bdefnn}[2]{\begin{defn}{\bf (#2)}\label{def:#1}\rm}
\newcommand{\brem}[1]{\begin{remark}\label{rem:#1}\rm}
\newcommand{\bexam}[1]{\begin{example}\label{ex:#1}\rm}
\newcommand{\bass}[1]{\begin{assumption}\label{ass:#1}\rm}
\newcommand{\ethm}{\end{thm}}
\newcommand{\eprop}{\end{prop}}
\newcommand{\elem}{\end{lem}}
\newcommand{\ecor}{\end{cor}}
\newcommand{\edefn}{\end{defn}}
\newcommand{\erem}{\end{remark}}
\newcommand{\eexam}{\end{example}}
\newcommand{\eass}{\end{assumption}}
\newcommand{\bpf}{\Proof}

\newcommand{\epf}{\eop}

\newcommand{\Thm}[1]{Theorem~\ref{th:#1}}
\newcommand{\Prop}[1]{Proposition~\ref{prop:#1}}
\newcommand{\Lem}[1]{Lemma~\ref{lem:#1}}
\newcommand{\Cor}[1]{Corollary~\ref{cor:#1}}

\newcommand{\Rem}[1]{Remark~\ref{rem:#1}}
\newcommand{\Exam}[1]{Example~\ref{ex:#1}}
\newcommand{\Ass}[1]{Assumption~\ref{ass:#1}}
\newcommand{\eq}[1]{\eqref{eq:#1}}

\newcommand{\E}{\mathbb{E}}
\newcommand{\R}{\mathbb{R}}
\newcommand{\N}{\mathbb{N}}

\newcommand{\gm}{\gamma}
\newcommand{\dl}{\delta}

\newcommand{\lm}{\lambda}
\newcommand{\ro}{\rho}
\newcommand{\Gm}{\Gamma}
\newcommand{\Dl}{\Delta}

\newcommand{\abs}[1]{\left| #1 \right|}
\newcommand{\Norm}[2]{\left\| #1 \right\|_{#2}}

\newcommand{\abra}[1]{\left( #1 \right)}
\newcommand{\bbra}[1]{\left\{ #1 \right\}}
\newcommand{\cbra}[1]{\left[ #1 \right]}

\newcommand{\bdry}{\partial}
\newcommand{\pdel}[2]{\dfrac{\bdry #1}{\bdry #2}}

\newcommand{\wg}{\wedge}
\newcommand{\nab}{\nabla}
\newcommand{\e}{\mathrm{e}}

\newcommand{\cH}{\mathcal{H}}

\newcommand{\sP}{\mathscr{P}}

\DeclareMathOperator{\Span}{Span}

\title{
Duality on gradient estimates and Wasserstein controls 
}
\author{Kazumasa Kuwada\footnote{
Partially supported by the JSPS fellowship for research abroad
}
}
\date{}
\begin{document}
\maketitle
\raggedbottom 
\begin{abstract}
We establish a duality 
between $L^p$-Wasserstein control and $L^q$-gradient estimate 
in a general framework. 
Our result extends a known result 
for a heat flow on a Riemannian manifold. 
Especially, 
we can derive a Wasserstein control of a heat flow 
directly from the corresponding gradient estimate 
of the heat semigroup 
without using any other notion of lower curvature bound. 
By applying our result to 
a subelliptic heat flow on a Lie group, 
we obtain a coupling of heat distributions 
which carries a good control of their relative distance. 
\end{abstract}
\begin{list}{\textbf{Key words:}}
{\setlength{\labelwidth}{72pt}
\setlength{\leftmargin}{72pt}}
\item
Wasserstein distance, 
gradient estimate, 
subelliptic diffusion, 
Ricci curvature 
\end{list}
\begin{list}{\textbf{Mathematics Subject Classification (2000):}}
{
  \setlength{\labelwidth}{72pt}
  \setlength{\leftmargin}{72pt}
}
\item
49L15, 49N15, 22E30, 60J60
\end{list}

\section{Introduction} 
\label{sec:intro} 

There are several ways 
to formulate a quantitative estimate 
on rate of convergence to equilibrium. 
By means of functional inequalities, 
an $L^q$-gradient estimate for a heat semigroup $P_t$ 
\begin{equation} \label{eq:gradient0} 
| \nab P_t f | (x) 
\le 
\e^{-kt} 
P_t ( | \nab f |^q ) (x)^{1/q}
\end{equation}
has been known to be a very powerful tool. 
It implies several functional inequalities 
such as Poincar\'{e} inequalities (when $q=2$) 
and logarithmic Sobolev inequalities (when $q=1$), 
which quantify convergence rates 
(see \cite{LSI2000,Bak97,BBBC08,Led_geom-Markov} 
and references therein). 
As a different approach to this problem, 
F.~Otto 
\cite{Otto_CPDE01}
discussed 
a contraction of 
$L^p$-Wasserstein distance 
\begin{equation} \label{eq:contraction0}
d^{\, W}_p ( \mu_t , \nu_t ) 
\le 
\e^{-kt} d^{\, W}_p ( \mu_0 , \nu_0 ) 
\end{equation}
for two (linear or nonlinear) diffusions 
$\mu_t , \nu_t$ of masses when $p=2$. 
His heuristic observation based on 
the geometry of the $L^2$-Wasserstein space 
has been a source of enormous developments 
in the theory of optimal transport
(see \cite{book_Vil2} and references therein). 
To investigate a relation 
between these formulations 
makes a connection between different approaches 
and hence it is an interesting problem. 
M.-K.~von Renesse and K.-Th.~Sturm 
\cite{Stu_Renes05} 
unified several formulations of this kind 
for linear heat equation 
on a complete Riemannian manifold. 
As a consequence of their work, 
\eq{gradient0} or \eq{contraction0} 
is shown to be equivalent to 
the presence of 
a lower Ricci curvature bound by $k$ 
(it also holds for $k<0$). 
But, in a more general framework, 
such a sort of duality has been known 
only when $p=1$ and $q=\infty$, 
which is the weakest form 
for \eq{gradient0} and \eq{contraction0} both. 

The main result of this paper 
extends the duality 
to that 
between an $L^q$-gradient estimate 
and 
an $L^p$-Wasserstein control 
for $p,q \in [1, \infty]$ 
with $p^{-1} + q^{-1} = 1$ 
beyond the case of 
a heat flow on a complete Riemannian manifold 
(see \Thm{duality} for the precise statement). 
We should emphasize that our duality does not require 
any other kind of curvature conditions. 
An $L^\infty$-Wasserstein control 
has been used in the literature 
as a tool to show $L^1$-gradient estimate 
in a coupling method 
for stochastic processes 
(for instance, see \cite{Wang97} 
and references therein). 
In the case of heat flows 
in a complete Riemannian manifolds, 
any construction of a coupling 
which carries $L^\infty$-Wasserstein control 
relies on lower Ricci curvature bounds. 
In fact, such an argument was used 
in von Renesse and Sturm's work. 
As a result,  
their proof employs 
a lower Ricci curvature bound 
to deduce 
Wasserstein controls 
from gradient estimates. 
Our result enables us to derive Wasserstein 
controls directly from gradient estimates. 
Such an implication is not known 
even in the case of heat flows on a Riemannian manifold. 
Furthermore, this is a great advantage under the lack of 
an appropriate notion of lower curvature bounds. 

Our work is strongly motivated 
by recent development on gradient estimates 
on a Lie group endowed 
with a sub-Riemannian structure 
\cite{BBBC08,BauBonn08,Dri-Mel05,Eldr09,HQLi06,Mel08}. 
To explain a consequence of our duality, 
we deal with 
the 3-dimensional Heisenberg group here. 
It is the simplest example of 
spaces possessing a non-Riemannian 
sub-Riemannian structure 
like a flat Euclidean space in Riemannian geometry. 
But, unlike Euclidean spaces, 
some results 
\cite{Dri-Mel05,Juil_CD-Heis} 
indicate 
that the ``Ricci curvature'' 
should be regarded as being unbounded from below 
(in a generalized sense). 
Nevertheless, $L^q$-gradient estimates 
hold for $q \in [ 1, \infty ]$ 
with a constant $K > 1$ 
instead of $\e^{-kt}$ in \eq{gradient0} 
\cite{BBBC08,Dri-Mel05,Eldr09,HQLi06}. 
We can apply our duality to this case 
to obtain the corresponding 
$L^p$-Wasserstein control 
for any $p \in [ 1 , \infty]$. 
In the theory of optimal transport
on the Heisenberg group, 
an $L^2$-Wasserstein control for the heat flow 
would be important 
(cf. \cite{Juil_gradF-Heis}). 
In probabilistic point of view, 
the heat flow is described 
by motions of a pair of 
the 2-dimensional Euclidean Brownian motion 
and 
the associated L\'{e}vy stochastic area. 
Our $L^\infty$-Wasserstein control means 
the existence of a coupling of 
two particles so that 
the distance between them at time $t$ is 
controlled by the initial distance almost surely. 
It is sometimes a complicated issue 
to construct a ``well-behaved'' coupling 
in the absence of curvature bounds. 
Especially, 
see \cite{BA-Crans-Kend95,Kend_Levy-couple} 
for works on a successful coupling 
on the Heisenberg group 
and its extension. 
Note that 
our formulation also fits with studying 
a heat semigroup under backward (super-)Ricci flow, 
in which case Wasserstein contractions 
with respect to a time-dependent distance function 
is shown recently 
\cite{Arn-Coul-Thal_horiz,McC-Topp_Wass-RF}. 

The notion of lower Ricci curvature bound 
has been extended in many ways. 
Although our result does not need those notions, 
they should be related 
since \eq{gradient0} and \eq{contraction0} 
are analytic and probabilistic characterizations 
of a lower Ricci curvature bound respectively. 
Here we review two extensions 
and observe how these are connected with our result. 
In an analytic way, D.~Bakry and M.~Emery 
\cite{Bak-Eme_CRAS84} 
(see also \cite{LSI2000} and references therein) 
extend the notion of lower Ricci curvature bound 
to $\Gm_2$-criterion or 
curvature-dimension condition. 
In an abstract framework where it works, 
a $\Gm_2$-criterion is equivalent to 
an $L^1$-gradient estimate. 
Note that their notion of gradient is 
different from ours. 
But, once these two notions coincide, 
a $\Gm_2$-criterion becomes equivalent to 
$L^\infty$-Wasserstein control with the aid of our result. 
In a sufficiently regular case as diffusions on a manifold, 
such an equivalence is well-known. 
Our result possibly provides an extension of this equivalence. 
In connection with the theory of optimal transport, 
convexities of entropy functionals are proposed 
by J.~Lott, C.~Villani and K.-Th.~Sturm 
\cite{Lott-Vill_AnnMath09,Sturm_Ric} 
as a natural extension of lower Ricci curvature bound. 
Under this condition, 
the existence 
of a heat flow 
and an $L^2$-Wasserstein control 
follow in some cases 
beyond Riemannian manifolds 
\cite{Ohta_grad-Alex,Sav_CRAS07} 
(see \cite{Erb_gradF-Riem,book_Vil2} 
for the case on a Riemannian manifold). 
With the aid of Theorem 8 in \cite{Sav_CRAS07}, 
we can apply our duality 
to show an $L^2$-gradient estimate 
for the heat semigroup. 
%
%

The idea of the proof of our main theorem is simple. 
The implication from 
a Wasserstein control to 
the corresponding gradient estimate 
is just a slight modification of existing arguments. 
The converse is based 
on the Kantorovich duality. 
If $p=1$, the Kantorovich duality becomes 
the Kantorovich-Rubinstein formula and 
the problem becomes much simpler. 
In the case $p > 1$, 
we employ a general theory of Hamilton-Jacobi 
semigroup developed in 
\cite{BEHM_HJ,Lott-Vill_HJ} 
to analyze the variational formula. 
When $p = \infty$, 
we use an approximation of $p$ by finite numbers 
because 
we are no longer able to apply 
the Kantorovich duality directly. 
Note that 
no semigroup property for heat semigroups 
is required in the proof. 
With keeping such a generality, 
our duality is sufficiently sharp in the sense 
that the control rate does not change 
when we obtain one estimate from the other, 
like the same $\e^{-kt}$ appears 
in \eq{gradient0} and \eq{contraction0} both. 

The organization of this paper is as follows. 
In the next section, we introduce 
our framework and state our main theorem. 
We review the notion of 
Wasserstein distance and gradient there. 
Our main theorem is shown in section \ref{sec:duality}. 
For the proof, we show basic properties 
of Wasserstein distances 
and summerize recent results 
on Hamilton-Jacobi semigroup there. 
In section \ref{sec:application}, 
we consider a heat flow on a sub-Riemannian manifold 
and apply our main theorem to these cases. 

\section{Framework and the main result}
\label{sec:pre}

Let $(X,d)$ be a complete, separable, proper, length metric space.  
Here, we say that $d$ is a length metric 
if, for every $x,y \in X$, 
$d(x,y)$ equals infimum of 
the length of a curve joining $x$ and $y$.  
Properness means that 
all closed metric balls in $X$ of finite radii are compact. 
Under these assumptions, 
there exists a curve joining $x$ and $y$ 
whose length realizes 
$d(x,y)$ for each $x,y$ (see \cite{BBI}, for instance). 
We call it minimal geodesic. 
%
Let $\tilde{d}$ be a continuous distance function on $X$, 
possibly different from $d$. 
Assume that for any $x,y \in X$, 
there is a minimal geodesic 
with respect to $\tilde{d}$ joining $x$ and $y$. 
We call such a curve ``$\tilde{d}$-minimal geodesic''. 
%
%

For two probability measures $\mu$ and $\nu$ on $X$, 
we denote the space of all couplings of $\mu$ and $\nu$ 
by $\Pi ( \mu , \nu )$. 
That is, 
$\pi \in \Pi ( \mu , \nu )$ means that 
$\pi$ is a probability measure on $X \times X$ 
satisfying 
$\pi ( A \times X ) = \mu (A)$ 
and  
$\pi ( X \times A ) = \nu (A)$ 
for each Borel set $A$. 
For $p\in [1, \infty]$ and a measurable function 
$\ro \: : \: X \times X \to [ 0, \infty)$, 
we define $\ro^{\, W}_p ( \mu , \nu )$ by 
\begin{equation} \label{eq:Wasserstein} 
\ro^{\, W}_p ( \mu , \nu ) 
: = 
\inf 
\bbra{ 
  \left. 
  \Norm{\ro}{L^p (\pi)}
  \; \right| \; 
  \pi \in \Pi ( \mu , \nu ) 
} . 
\end{equation} 
We are interested 
in the case $\ro = d$ and $\ro = \tilde{d}$. 
If $d^{\, W}_p ( \mu, \nu ) < \infty$, 
then there always exists 
a minimizer of the infimum 
on the right hand side in \eq{Wasserstein}. 
In addition, 
$d^{\, W}_p$ satisfies all properties of distance function 
on the space of probability measures 
though it may take the value $+\infty$. 
The same are also true for $\tilde{d}^{\, W}_p$. 
These facts are well-known for $p \in [1, \infty)$ 
and 
we can show it similarly even when $p = \infty$. 
It is sometimes reasonable to restrict $d^{\, W}_p$ 
on all probability measures having finite $p$-th moments 
in order to ensure $d^{\, W}_p ( \mu , \nu ) < \infty$. 
But, in this paper, we do not adopt such a restriction. 
Note that, when $p < \infty$, 
we usually call the restriction of $d^{\, W}_p$ 
the $L^p$-Wasserstein distance. 
See \cite{book_Vil1} 
for more details and a proof of these facts. 
%

Let $C_b (X)$ be 
the space of bounded continuous functions on $X$ 
equipped with the supremum norm. 
Let $C_{L} (X)$ be the collection of 
all Lipschitz continuous functions on $X$ and 
$C_{b,L} (X) := C_b (X) \cap C_L (X)$. 
Note that, if we merely say ``Lipschitz'', 
it means ``Lipschitz with respect to $d$''. 
For Lipschitz continuity 
with respect to $\tilde{d}$, 
we use the expression ``$\tilde{d}$-Lipschitz''. 

For a measurable function $f$ on $X$ 
and $x \in X$, 
we define $| \nab_{d} f | (x)$ by 
\[
| \nab_{d} f | (x) 
= 
\lim_{r \downarrow 0} 
\sup_{0 < d(x,y) \le r}
\abs{ 
  \frac{f (x) - f (y)}{d(x,y)} 
}
.
\] 
We set 
$\| \nab_{d} f \|_\infty = \sup_{x \in X} | \nab_{d} f |(x)$.   
Note that 
$\| \nab_{d} f \|_\infty < \infty$ holds 
if and only if $f \in C_L (X)$. 
In addition, for $f \in C_L (X)$, 
\begin{equation} \label{eq:Lipschitz-bound}
\| \nab_{d} f \|_\infty 
= 
\sup_{x \neq y} 
\abs{ \frac{ f(x) - f(y) }{ d(x,y) } }
. 
\end{equation} 
For a pair of measurable functions $f$ and $g$ on $X$, 
we say that $g$ is an upper gradient of $f$ 
if, for each rectifiable curve 
$\gm \: : \: [ 0, l ] \to X$ 
parametrized with the arc-length, we have 
\[
\abs{ 
  f ( \gm (l) ) - f ( \gm (0) ) 
}
\le 
\int_0^l g ( \gm (s) ) ds 
.
\] 
We will use the following fact as a basic tool. 
\blemn{upper-gradient}
{\protect{\cite[Proposition~1.11]{Chee99},
\cite[Proposition~10.2]{SmP}}} 
For $f \in C_L (X)$, 
$| \nab_{d} f |$ is an upper gradient of $f$. 
\elem  
We also use the same notations for $\tilde{d}$. 
All the properties described above for $| \nab_{d} f |$, 
including \Lem{upper-gradient}, 
are also true for $| \nab_{\tilde{d}} f |$. 

Set $\sP (X)$ be the space of 
all probability measures on $X$ 
equipped with the topology of weak convergence. 
Let $( P_x )_{x \in X}$ be a family of elements in $\sP (X)$. 
Assume that $x \mapsto P_x$ is continuous 
as a map from $X$ to $\sP (X)$. 
Then $( P_x )_{x \in X}$ defines 
a bounded linear operator $P$ on $C_b (X)$ 
by $Pf (x) : = \int_X f(y) P_x (dy)$. 
Let $P^*$ be the adjoint operator of $P$. 
Note that $P^* ( \sP (X) ) \subset \sP (X)$ holds. 

For describing our main theorem, 
we state the following conditions: 
\bass{regularity} 
There exists a positive Radon measure $v$ on $X$ such that 
\begin{enumerate}
\item \label{doubling}
$(X,d,v)$ enjoys the local volume doubling condition. 
That is, there are constants $D , R_1 > 0$ 
such that 
$v ( B_{2r} ( x ) ) \le D v ( B_{r} (x) )$ holds  
for all $x \in X$ and $r \in ( 0 , R_1 )$. 
\item \label{Poincare}
$(X,d,v)$ supports a $(1,p_0)$-local Poincar\'{e} inequality 
for some $p_0 \ge 1$. That is, 
for every $R > 0$, 
there are constants 
$\lm \ge 1$ and $C_P > 0$ 
such that, 
for any $f \in L^1_{\mathrm{loc}} (v)$ and 
any upper gradient $g$ of $f$, 
\begin{equation} \label{eq:Poincare}
\int_{B_r (x)} \abs{ f - f_{x,r} } d v 
\le 
C_P r 
\bbra{ 
  \int_{ B_{\lm r} (x)} g^{p_0}  dv 
}^{1/{p_0}}
\end{equation}
holds 
for every $x \in X$ and $r \in ( 0 , R )$,  
where 
$f_{x,r} := v ( B_r (x) )^{-1} \int_{B_r (x)} f \; d v$. 
\item \label{density}
$P_x$ is absolutely continuous 
with respect to $v$ for all $x \in X$; 
$P_x (dy) = P_x (y) v (dy)$. 
In addition, 
the density $P_x(y)$ is continuous 
with respect to $x$. 
\end{enumerate} 
\eass
Now we are in turn to state our main theorem. 
\bthm{duality} 
Suppose that \Ass{regularity} holds. 
Then, 
for any $p \in [1,\infty]$, 
the following are equivalent; 
\begin{enumerate}
\item \label{transport-bound}
For all $\mu , \nu \in \sP (X)$, 
\begin{equation} 
\tag{$C_p$}
\label{eq:cost} 
d^{\, W}_p ( P^* \mu , P^* \nu  ) 
\le 
\tilde{d}^{\, W}_p ( \mu , \nu ) 
. 
\end{equation}
\item \label{gradient-estimate}
When $p > 1$, for all $f \in C_{b,L} (X)$ and $x \in X$, 
\begin{equation} 
\tag{$G_q$}
\label{eq:gradient}
| \nab_{\tilde{d}} P f | (x) 
\le 
P ( | \nab_{d} f |^q )(x)^{1/q} 
, 
\end{equation}
where $q$ is the H\"{o}lder conjugate of $p$; $1/p + 1/q =1$.  
When $p=1$, for all $f \in C_{b,L} (X)$, 
\begin{equation} 
\tag{$G_\infty$}
\label{eq:gradient-infinity}
\Norm{ \nab_{\tilde{d}} P f }{\infty} 
\le 
\Norm{ \nab_{d} f }{\infty} 
. 
\end{equation} 
\end{enumerate} 
\ethm 
\brem{Poincare} 
We give several remarks on \Ass{regularity} and \Thm{duality}. 
\begin{enumerate}
\item 
If \Ass{regularity} \ref{doubling} holds, 
then \Ass{regularity} \ref{Poincare} follows 
once we obtain \eq{Poincare} with $p_0 = 1$ 
for \emph{some} $R >0$ 
by a well-known argument. 
See \cite[Lemma~5.3.1]{book_Sal}, 
for instance. 
The same is true for a (2,2)-Poincar\'{e} inequality, 
which yield a (1,2)-Poincar\'{e} inequality. 
\item
It is shown in \cite{Chee99} that, 
under 
\Ass{regularity} \ref{doubling} \ref{Poincare}, 
$| \nab_{d} f |$ coincides with 
an $L^{p_0}$-minimal generalized upper gradient $g_f$ 
for those $f$ for which $g_f$ is well-defined. 
This fact itself is not used in this article. 
But, it will be helpful 
when we apply our main theorem 
to more concrete problems. 
In fact, 
the notion of minimal generalized upper gradients 
is regarded as a sort of weak derivative 
in the theory of Sobolev spaces. 
We can identify these two notions 
on Euclidean spaces or Riemannian manifolds. 
\item
\Ass{regularity} is used only when we show the implication 
\eq{gradient} $\Rightarrow$ \eq{cost}
for $p \in ( 1 , \infty ]$. 
Thus the rest holds true without \Ass{regularity}. 
We need \Ass{regularity} \ref{doubling} \ref{Poincare} 
only for employing a property of 
Hamilton-Jacobi semigroups. 
To make these facts clear, in the rest of this paper, 
we will mention \Ass{regularity} when we require it. 
\item
The duality between \eq{gradient0} and \eq{contraction0} 
is resumed by choosing $P = P_t$ and $\tilde{d} = \e^{-kt} d$. 
The case $\tilde{d}$ is essentially different from $d$ 
naturally occurs 
if we consider 
a heat flow under a backward (super-)Ricci flow 
(see \cite{Arn-Coul-Thal_horiz,McC-Topp_Wass-RF}).
\item 
Obviously ($G_q$) implies ($G_{q'}$) 
for $q, q' \in [ 1 , \infty ]$ with $q < q'$ 
by the H\"{o}lder inequality.  
The dual implication ($C_p$) $\Rightarrow$ ($C_{p'}$) 
for $p, p' \in [ 1, \infty ]$ with $p > p'$ 
also holds true 
without using the equivalence in \Thm{duality} 
(see \Cor{monotone} below). 
For a heat flow on a Riemannian manifold 
(i.e.~$P = P_t$ and $\tilde{d}= \e^{-kt}$), 
if \eq{cost} or \eq{gradient} holds 
for \emph{some} $p \in [1,\infty]$, 
then \eq{cost} and \eq{gradient} hold 
for \emph{any} $p \in [1,\infty]$. 
At this moment, 
it is not clear that what condition 
guarantees such a ``$L^p$-independence''. 
\end{enumerate}
\erem
\section{Proof of \protect{\Thm{duality}}}
\label{sec:duality} 

We begin with 
showing 
the implication \eq{cost} $\Rightarrow$ \eq{gradient}. 
\bprop{cost-gradient} 
Suppose \eq{cost} for $p \in [1, \infty]$. 
Then \eq{gradient} holds for $q \in [1,\infty]$ 
with $p^{-1} + q^{-1} = 1$. 
\eprop 
\bpf 
For $x,y \in X$, 
take $\pi_{xy} \in \Pi ( P_x , P_y )$ 
such that 
$\| d \|_{L^p (\pi_{xy})}  =  d^{\, W}_p (P_x , P_y )$. 
Since $P_z = P^* \dl_z$ for $z \in X$, 
\eq{cost} yields 
$
d^{\, W}_p ( P_x , P_y ) 
\le 
\tilde{d}^{\, W}_p ( \dl_x , \dl_y ) 
= \tilde{d} (x,y) 
$. 
For $f \in C_{b, L} (X)$, 
\begin{align*}
\abs{ 
  P f (x) - P f (y) 
} 
& = 
\abs{ 
  \int_X f \, d P_x  
  - 
  \int_X f \, d P_y 
} 
\le 
\int_{X \times X} 
\abs{ 
  f (z) - f (w) 
}
\pi_{xy} ( dz dw ). 
\end{align*} 

\textbf{(i) The case $p=1$:} 
\eq{Lipschitz-bound} 
together with ($C_1$)
implies 
\[
\int_{X \times X} 
\abs{ 
  f (z) - f (w) 
}
\pi_{xy} ( dz dw ) 
\le 
\Norm{ \nab_{d} f }{\infty}
d^{\, W}_1 ( P_x , P_y ) 
\le 
\Norm{ \nab_{d} f }{\infty}
\tilde{d} (x,y) . 
\]
Hence, by dividing the above inequalities by $\tilde{d}(x,y)$ and 
by taking supremum in $x \neq y$, 
the conclusion follows. 

\textbf{(ii) The case $p \in (1, \infty)$:} 
Let us define $G_r \: : \: X \to \R$ 
by 
\begin{align*}
G_r (z) : = 
\sup_{w \in B_r (z) \setminus \{ z \} } 
\abs{ 
  \frac{ f (z) - f (w) }{ d (z,w) } 
}. 
\end{align*}
Set $r := \tilde{d} (x,y)^{1/(2q)}$. 
The H\"{o}lder inequality 
and 
the Chebyshev inequality 
yield 
\begin{align*} 
\int_{X \times X} 
& | f (z) - f (w) | 
\pi_{xy} ( dz dw ) 
\\
& = 
\int_{X \times X} 
\abs{ 
  \frac{f (z) - f (w)}{d (z,w)} 
} 
1_{ \{ 0 < d (z,w) \le r \} }
d (z,w)
\pi_{xy} ( dz dw ) 
\\
& \qquad + 
\int_{X \times X} 
\abs{ 
  f (z) - f (w)
}  
1_{ \{ d (z,w) > r \} }
\pi_{xy} ( dz dw ) 
\\
& \le 
\bbra{ 
  \int_{X \times X} 
  \abs{ 
    \frac{ f(z) - f(w) }{ d (z,w) } 
  }^q 
  1_{ \{ 0 < d (z,w) \le r \} }
  \pi_{xy} ( dz dw ) 
}^{1/q}
\| d \|_{L^p ( \pi_{xy} )} 
+ 
\frac{ 2 \| f \|_{\infty} \| d \|_{L^p (\pi_{xy})}^p }{ r^p }
\\ 
& \le 
\Norm{ G_r }{L^q ( P_x )}
d^{\, W}_p ( P_x , P_y ) 
+ 
\frac{ 2 \| f \|_{\infty} d^{\, W}_p ( P_x ,P_y )^p }{ r^p }
\\
& \le 
\Norm{ G_r }{L^q ( P_x )}
\tilde{d} (x,y)
+ 
2 \| f \|_{\infty} \tilde{d} ( x ,y )^{1+(p-1)/2} 
. 
\end{align*}
Here the last inequality follows from \eq{cost}. 
Since $\lim_{y \to x} r = 0$, 
$\lim_{y \to x} G_r(z) = \abs{ \nab_{d} f } (z)$ holds. 
By virtue of 
$\abs{ G_r (z) } \le \Norm{ \nab_{d} f }{ \infty }$, 
we can apply the dominated convergence theorem 
to obtain 
$
\lim_{y \to x} \Norm{G_r}{L^q (P_x)} 
= 
\Norm{ \abs{ \nab_{d} f } }{L^q ( P_x )} 
$. 
Thus, 
by dividing the above inequalities by $\tilde{d} (x,y)$ and 
by tending $y \to x$, the conclusion follows. 

\textbf{(iii) The case $p = \infty$:} 
($C_\infty$) implies $d (z,w) \le \tilde{d} (x,y)$ 
for $\pi_{xy}$-a.e.~$(z,w)$. 
Hence we have 
\[
\int_{X \times X} | f(z) - f(w) | \pi_{xy} ( dz dw ) 
\le 
\tilde{d} (x,y) \| G_{\tilde{d}(x,y)} \|_{ L^1 ( P_x ) }
. 
\]  
Thus the proof will be completed 
by following a similar argument as above. 
\epf

For the converse implication, 
first we show two auxiliary lemmas 
concerning to Wasserstein distances. 
The first one will be used 
to deal with $L^\infty$-Wasserstein distance. 
\blem{infinite-cost}
Let $\ro \: : \: X \times X \to [0,\infty)$ be 
a continuous function. 
Then 
$
\lim_{p \to \infty} \ro^{\, W}_p ( \mu , \nu ) 
= 
\ro^{\, W}_\infty ( \mu , \nu ) 
$ 
for any $\mu , \nu \in \sP (X)$. 
\elem
\bpf
Note that $\ro^{\, W}_p ( \mu , \nu )$ is increasing in $p$ 
by the H\"{o}lder inequality. 
Hence 
$C := \lim_{p \to \infty} \ro^{\, W}_p (\mu , \nu ) \in [0, \infty]$ 
exists. 
Take 
$\pi_n \in \Pi (\mu, \nu)$ for $n \in \N$ 
such that 
$\ro^{\, W}_{n} ( \mu , \nu ) = \Norm{ \ro }{L^{n} ( \pi_n )}$ 
hold. 
Since $\pi_n \in \Pi ( \mu , \nu )$, 
$( \pi_n )_{n \in \N}$ is tight. 
Thus there exists a convergent subsequence 
$( \pi_{n_k} )_{k \in \N}$ of 
$( \pi_n )_{n \in \N}$. 
We denote the limit of $\pi_{n_k}$ by $\pi_\infty$. 
Take $R >0$ and $n \in \N$ arbitrary. 
Since $\ro \wg R \in C_b (X \times X)$, 
we have 
\[
\Norm{\ro \wg R}{L^{n} ( \pi_\infty )} 
= 
\lim_{k \to \infty} \Norm{\ro \wg R}{L^{n} ( \pi_{n_k} )} 
\le 
\lim_{k \to \infty} \Norm{ \ro }{L^{n_k} ( \pi_{n_k} )} 
= 
C . 
\]
Here the inequality follows from 
the H\"{o}lder inequality for sufficiently large $k$. 
Thus, as $R \to \infty$ and $n \to \infty$, 
we obtain 
$\Norm{\ro}{L^\infty ( \pi_\infty )} \le C$. 
Thus the assertion holds 
if $\ro^{\, W}_\infty ( \mu , \nu ) = \infty$.  
When $\ro^{\, W}_\infty ( \mu , \nu ) < \infty$, 
we can take $\pi \in \Pi ( \mu , \nu )$ 
such that $\Norm{ \ro }{ L^\infty ( \pi ) } < \infty$.   
Then 
$\ro^{\, W}_p ( \mu , \nu ) 
\le \Norm{ \ro }{L^p ( \pi )} 
\le \Norm{ \ro }{ L^\infty ( \pi ) } 
$. 
Thus $C \le \Norm{ \ro }{ L^\infty ( \pi )}$ holds. 
It yields $C \le \ro^{\, W}_\infty ( \mu, \nu )$ and 
hence the conclusion holds. 
\epf
The next one is useful to reduce the problem 
in a simpler case. 
\blem{From-Dirac} 
If \eq{cost} holds 
for any pair of Dirac measures, 
then \eq{cost} holds 
for any $\mu , \nu \in \sP (X)$.  
\elem
Although this is probably well-known for experts 
at least when $p \in [1, \infty)$, 
we give a proof for completeness. 
\bpf
First we consider the case $p < \infty$. 
Given $\mu , \nu \in \sP (X)$, 
take $\pi \in \Pi ( \mu , \nu )$ 
so that 
$
\| \tilde{d} \|_{L^p (\pi)} 
= 
\tilde{d}^{\, W}_p ( \mu, \nu )
$. 
We may assume $\tilde{d}^{\, W}_p (\mu , \nu ) < \infty$ 
without loss of generality. 
For $x,y \in X$, 
take $P_{x,y} \in \Pi ( P_x , P_y )$ 
so that 
$
\| d \|_{ L^p ( P_{x,y} )} 
= 
d^{\, W}_p ( P_x , P_y )
$. 
By Corollary 5.22 of \cite{book_Vil2}, 
we can choose $\{ P_{x,y} \}_{x,y \in X}$ 
so that 
the map $(x,y) \mapsto P_{x,y}$ is measurable. 
Define $\tilde{\pi} \in \Pi ( P^* \mu , P^* \nu )$ 
by $\tilde{\pi}(A) := \int_{X \times X} P_{x,y} (A) \pi (dx dy)$. 
Then \eq{cost} for Dirac measures implies   
\begin{align*}
d^{\, W}_p ( P^* \mu , P^* \nu ) 
& \le 
\| d \|_{L^p ( \tilde{\pi} ) }
= 
\bbra{ 
  \int_{X \times X} 
  \| d \|_{ L^p ( P_{x,y} ) }^p 
  \pi ( dx dy ) 
}^{1/p} 
\le 
\| \tilde{d} \|_{ L^p (\pi) } 
= 
\tilde{d}^{\, W}_p ( \mu , \nu ). 
\end{align*} 
Thus the assertion holds. 
When $p = \infty$, ($C_\infty$) 
for Dirac measures implies 
($C_{p'}$) for Dirac measures 
for any $1 \le p' < \infty$. 
Thus we obtain ($C_{p'}$) 
for any $\mu ,\nu \in \sP (X)$. 
Hence applying \Lem{infinite-cost} 
for $\ro = d$ and $\ro = \tilde{d}$ 
yields 
the conclusion. 
\epf
By the H\"{o}lder inequality, 
$\eq{cost}$ for Dirac measures  
yields ($C_{p'}$) for Dirac measures 
if $p' < p$. 
Thus we obtain the following as a by-product of \Lem{From-Dirac}.  
\bcor{monotone}
\eq{cost} implies $(C_{p'})$ 
for any $p, p' \in [ 1, \infty]$ with $p > p'$. 
\ecor

Next we introduce 
the notion and some properties of 
Hamilton-Jacobi semigroup, 
which plays an essential role in the sequel. 
Let $L \: : \: [0,\infty) \to [ 0, \infty)$ be 
a convex superlinear function with $L(0)=0$.   
Note that $L$ is continuous and increasing.  
We denote the Legendre conjugate of $L$ 
by $L^* \: : \: [0,\infty) \to [ 0, \infty)$, 
which is given by 
$L^* (z) = \sup_{w \ge 0} \cbra{ wz - L (w)}$. 
For $f \in C_b (X)$ and $t >0$, 
we define a function $Q_t f$ on $X$ 
by 
\[
Q_t f (x) 
: = 
\inf_{y \in X}  
\cbra{ 
  f(y) 
  + 
  t L 
  \abra{ 
    \frac{d(x,y)}{t}
  } 
} 
.
\] 
For convenience, we write $Q_0 f := f$. 
We call $Q_t$ 
the Hamilton-Jacobi semigroup associated with $L$.  
%
Several basic properties of $Q_t f$ 
in an abstract framework 
are studied 
in \cite{BEHM_HJ,Lott-Vill_HJ}. 
In \cite{Lott-Vill_HJ}, 
they assumed $X$ to be compact and $L(s) = s^2$. 
In \cite{BEHM_HJ}, they assumed $f \in C_L (X)$. 
Among them, 
the following are all we need in this paper. 
\blemn{Hamilton-Jacobi}
{\protect{\cite[Theorem~2.5]{BEHM_HJ}, 
\cite[Theorem~2.5]{Lott-Vill_HJ}}}
\begin{enumerate}
\item \label{HJ-bound}
$\inf_{y \in X} f (y) \le Q_t f (x) \le f (x)$.  
In particular, $Q_t f \in C_b (X)$. 
\item \label{semigroup}
$Q_t ( Q_s f ) = Q_{t+s} f$. 
\item \label{HJ-initial}
$Q_t f(x)$ is nonincreasing in $t$ and 
$\lim_{t \downarrow 0} Q_t f (x) = f (x)$. 
\item \label{HJ-Lipschitz}
Set $u (t,x) = Q_t f(x)$. 
If $f \in C_L (X)$, 
then $u \in C_L ((0,\infty) \times X)$.  
Moreover, 
\[
\sup_{
  \begin{subarray}{c} 
    s \neq t 
    \\ 
    y \neq x 
  \end{subarray}
} 
\frac{\abs{ u (t,x) - u (s,y) }}{|t-s| + d(x,y)} 
\le 
\Norm{ \nab_{d} f }{\infty}
\vee
L^* ( \Norm{ \nab_{d} f }{\infty} )
.
\] 
\item \label{Hamilton-Jacobi}
Suppose \Ass{regularity} 
\ref{doubling} \ref{Poincare}. 
Then, for $t > 0$ and $v$-a.e. $x \in X$, 
$Q_t f$ satisfies 
the Hamilton-Jacobi equation 
associated with $L^*$: 
\[
\lim_{s \downarrow 0} 
\frac{Q_{t+s} f (x) - Q_t f (x)}{s} 
= 
L^* ( \abs{ \nab_d Q_t f } (x) ) . 
\]
\end{enumerate}
\elem
We do not use 
\Lem{Hamilton-Jacobi} \ref{HJ-bound} \ref{semigroup} 
in the sequel. 
But, 
it explains why we call $Q_t$ ``semigroup'' well. 
Note that 
\Lem{Hamilton-Jacobi} \ref{Hamilton-Jacobi} is shown 
in \cite{BEHM_HJ,Lott-Vill_HJ} 
for the subgradient norm instead of 
the gradient norm $\abs{\nab_d f}$. 
Since these two notions coincides 
$v$-almost everywhere in this case 
(see \cite[Remark~2.27]{Lott-Vill_HJ}), 
\Lem{Hamilton-Jacobi} \ref{Hamilton-Jacobi} 
is still valid. 
 
Finally, we review 
the Kantorovich duality 
(see 
\cite[Theorem~1.3]{book_Vil1} or 
\cite[Theorem~5.10]{book_Vil2}, 
for example). 
For $\mu, \nu \in X$ and $1 \le p < \infty$, 
the following duality holds: 
\begin{align} 
\nn 
d^{\, W}_p ( \mu , \nu )^p 
& =
\sup 
\bbra{ 
  \left.   
  \int_X g \, d \mu - \int_X f \, d \nu 
  \; \right| \; 
  \begin{array}{l} 
    f,g \in C_b (X), 
    \\
    g(y) - f(x) \le d (x,y)^p 
    \mbox{ for all $x, y \in X$}
  \end{array}
}, 
\\
\label{eq:Kantorovich}
& = 
\sup_{f \in C_b (X)} 
\cbra{ 
  \int_X f^* \, d \mu - \int_X f \, d \nu 
}, 
\end{align}
where $f^* (y) : = \inf_{x \in X} \cbra{ f (x) + d (x,y)^p }$. 
In particular, when $p=1$, 
\eq{Kantorovich} is written as follows: 
\begin{equation} \label{eq:Rubinstein}
d^{\, W}_1 ( \mu , \nu ) 
= 
\sup_{
  \begin{subarray}{c} 
    f \in C_L (X)
    \\ 
    \Norm{ \nab f }{\infty} \le 1 
  \end{subarray} 
} 
\cbra{
  \int_M f \, d \mu - \int_M f \, d \nu 
}. 
\end{equation} 
This is so-called the Kantorovich-Rubinstein formula 
(see \cite[Theorem~1.14]{book_Vil1} or 
\cite[Particular Case~5.16]{book_Vil2}). 
\brem{Lipschitz-approximation}
An observation on the proof in \cite{book_Vil2} 
tells us that 
the latter supremum in \eq{Kantorovich} 
can be approximated 
by elements in $C_{b,L} (X)$. 
Actually, in that proof, there appears 
a sequence of pair of functions 
$\phi_k ,\psi_k \in C_b (X)$ 
approximating the former supremum in \eq{Kantorovich} 
by taking $f = \psi_k , g = \phi_k$. 
We can easily verify $\psi_k \in C_{b,L} (X)$ 
and that $( \psi_k )_{k \in \N}$ also approximates 
the latter supremum in \eq{Kantorovich}. 
Moreover, 
we can assume that 
each element of approximating sequence 
has a compact support without loss of generality, 
thanks to the tightness of $\mu, \nu$ 
and the properness of $X$. 
\erem
%
%
%

Now we are in position 
to complete the proof of \Thm{duality}. 
\bprop{gradient-cost} 
Suppose that \Ass{regularity} holds. 
Then \eq{gradient} implies \eq{cost} 
for $p, q \in [ 1, \infty ]$ with $p^{-1} + q^{-1} = 1$.  
\eprop
\bpf
By virtue of \Lem{From-Dirac}, 
it suffices to show \eq{cost} 
for $\mu = \dl_x , \nu = \dl_y$, $x \neq y$.  
Take a $\tilde{d}$-minimal geodesic 
$\gm \: : \: [ 0, 1 ] \to X$ 
from $y$ to $x$, which is re-parametrized 
to have a constant speed. 
Here ``constant speed'' means 
$\tilde{d} ( \gm_s , \gm_t ) = |s-t| \tilde{d} (x,y)$. 
Note that, by \eq{gradient}, 
$Pf$ is $\tilde{d}$-Lipschitz continuous 
if $f \in C_L (X)$. 

\textbf{(i) The case $p=1$:} 
The Kantorovich-Rubinstein formula 
\eq{Rubinstein} yields 
\begin{align} 
\label{eq:Rubinstein2}
d^{\, W}_1 ( P_x , P_y ) 
= 
\sup_{ 
  \begin{subarray}{c}
    f \in C_L (X) 
    \\
    \Norm{ \nab_{d} f }{\infty} \le 1 
  \end{subarray} 
} 
  \cbra{ P f (x) - P f (y) }
. 
\end{align}
For $f \in C_L (X)$, 
we can apply \Lem{upper-gradient} to $Pf$. 
Thus ($G_\infty$) yields 
\[
\abs{ P f (x) - P f (y) } 
\le 
\int_0^{\tilde{d} (x,y)} 
\abs{ \nab_{\tilde{d}} P f } ( \gm_s ) ds 
\le 
\Norm{ \nab_{d} f }{ \infty } \tilde{d} (x,y)
. 
\]
Combining this estimate with \eq{Rubinstein2}, 
the conclusion follows. 

\textbf{(ii) The case $1 < p < \infty$:}
Let $Q_t$ be the Hamilton-Jacobi semigroup 
associated with $L(s) := p^{-1} s^p$. 
Note that 
its Legendre conjugate $L^*$ is 
computed as $L^* (s) = q^{-1} s^q$. 
By 
\eq{Kantorovich} 
and 
\Rem{Lipschitz-approximation}, 
we have 
\begin{align}
d^{\, W}_p ( P_x , P_y )^p
& = 
\sup_{f \in C_{b,L}(X)} 
\cbra{ 
  P ( f^* ) (x) - P f (y)  
} 
\label{eq:Kantorovich2}
=
p \sup_{f \in C_{b,L}(X)} 
\cbra{ 
  P Q_1 f (x) - P f (y) 
} 
. 
\end{align} 
To obtain an integral expression 
of the term in the above supremum 
(see \eq{fundamental} below), 
we give some estimates. 
\eq{gradient} and 
\Lem{Hamilton-Jacobi} \ref{HJ-Lipschitz} 
yield 
\[
\abs{ \nab_{\tilde{d}} P Q_s f } ( z ) 
\le 
\Norm{ \abs{ \nab_{d} Q_s f } }{ L^q ( P_z ) }  
\le 
\Norm{ \nab_{d} f }{\infty} 
\vee 
L^* ( \Norm{ \nab_{d} f }{\infty} )  
\]
for $s \ge 0$ and $z \in X$.  
Thus \Lem{upper-gradient} and 
\Lem{Hamilton-Jacobi} \ref{HJ-Lipschitz} 
imply 
\begin{align*}
\abs{ 
  \frac{
    P Q_{t+s} f ( \gm_{t+s} ) 
    - 
    P Q_s f ( \gm_s ) 
  }{t}
}
& \le 
\abs{
  \frac{
    P Q_{t+s} f ( \gm_{t+s} ) 
    - 
    P Q_{t+s} f ( \gm_s ) 
  }{t}
}
+ 
\abs{
  \int_X 
  \frac{
    Q_{t+s} f 
    - 
    Q_s f 
  }{t} 
  d P_{ \gm_s }
}
\\
& \le 
\frac{\tilde{d} (x,y)}{t} 
\int_s^{t+s}
   \abs{ \nab_{\tilde{d}} P Q_{t+s} f }( \gm_{u} ) 
d u
+ 
\int_X 
\abs{
  \frac{
    Q_{t+s} f 
    - 
    Q_s f 
  }{t} 
}
d P_{ \gm_s }
\\
& \le 
\abra{ 1 + \tilde{d} (x,y) } 
\abra{  
  \Norm{ \nab_{d} f }{\infty} 
  \vee 
  L^* ( \Norm{ \nab_{d} f }{\infty} ) 
}
\end{align*}
for $s \ge 0$. 
It means that 
$P Q_s f ( \gm_s )$ is 
Lipschitz continuous 
as a function of $s \in [0,1]$. 
Hence there exists a derivative 
$\partial_s ( P Q_s f ( \gm_s ) )$ 
for a.e.$s \in [0,1]$ 
and we have 
\begin{equation} \label{eq:fundamental}
P Q_1 f (x) - P f (y) 
= 
\int_0^1 
\partial_s 
\abra{ 
  P Q_s f ( \gm_s ) 
} ds 
. 
\end{equation} 
Let $s \in (0,1)$ be a point 
where $P Q_s f ( \gm_s )$ is differentiable. 
It implies 
\begin{align} \label{eq:derivative0}
\partial_s ( P Q_s f ( \gm_s ) ) 
& = 
\lim_{t \downarrow 0} 
\frac{ 
  P Q_{s+t} f ( \gm_{s+t} ) 
  - 
  P Q_s f ( \gm_s )
}{t} 
. 
\end{align}
We have  
\begin{align}
\frac{ 
  P Q_{s+t} f ( \gm_{s+t} ) 
  - 
  P Q_s f ( \gm_s )
}{t} 
& = 
\int_X 
\frac{ 
   Q_{s+t} f 
  - 
   Q_s f 
}{t}
d P_{ \gm_{s+t} } 
\label{eq:derivative}
+ 
\frac{ 
  P Q_s f ( \gm_{s+t} ) 
  - 
  P Q_s f ( \gm_s )
}{t}
. 
\end{align}
By 
\Lem{upper-gradient} 
together with \eq{gradient}, 
\begin{align} 
\frac{ 
  P Q_s f ( \gm_{s+t} ) 
  - 
  P Q_s f ( \gm_s )
}{t}
\label{eq:difference1}
& \le 
\frac{ \tilde{d} (x,y) }{t} 
\int_s^{s+t} 
\bbra{ 
  \abra{ 
    P ( \abs{ \nab_{d} Q_s f }^q 
  } ( \gm_u ) 
}^{1/q} 
du 
. 
\end{align} 
By virtue of \Ass{regularity} \ref{density}, 
the Fatou lemma together with 
the boundedness of $\abs{ \nab_{d} Q_t f }$ 
implies that 
$( P | \nab_d Q_s f |^q )( \gm_u )$ is 
upper semi-continuous in $u$. 
Thus \eq{difference1} yields 
\[
\limsup_{t \downarrow 0} 
\frac{ 
  P Q_s f ( \gm_{s+t} ) 
  - 
  P Q_s f ( \gm_s )
}{t}
\le 
\tilde{d} (x,y) 
\Norm{ \abs{ \nab_{d} Q_s f } }{L^q ( P_{\gm_s} )}
.   
\]
For the first term in \eq{derivative}, 
\Lem{Hamilton-Jacobi} \ref{HJ-initial} 
implies the integrand is nonpositive. 
Thanks to \Ass{regularity} \ref{doubling} \ref{Poincare}, 
\Lem{Hamilton-Jacobi} \ref{Hamilton-Jacobi} 
is applicable to the integrand. 
Thus the Fatou lemma together with 
\Ass{regularity} \ref{density} 
yields 
\begin{align}
\nn 
\limsup_{t \downarrow 0} 
\int_X 
\frac{ Q_{t+s} f - Q_s f }{t} dP_{\gm_{s+t}}
& = 
\limsup_{t \downarrow 0} 
\int_X 
\frac{ Q_{t+s} f(z) - Q_s f (z) }{t} 
P_{\gm_{s+t}} (z) v (dz) 
\\
\nn 
& \le 
\int_X 
\limsup_{t \downarrow 0} 
\frac{ Q_{t+s} f(z) - Q_s f (z) }{t} 
P_{\gm_{s+t}} (z) 
v (dz)
\\
\label{eq:difference2} 
& = 
- 
\int_X 
L^* 
\abra{ 
  \abs{ 
    \nab_{d} Q_s f
  } (z)
} 
P_{\gm_s} (z) 
v (dz)
.  
\end{align} 
Combining 
\eq{derivative}, 
\eq{difference1} 
and 
\eq{difference2} 
with \eq{fundamental} and \eq{derivative0}, 
\begin{align*}
\nn
P Q_1 f (x) - P f (y) 
& \le 
\int_0^1 
\abra{ 
  \tilde{d} (x,y) 
  \Norm{ 
    \abs{ 
      \nab_{d} Q_s f 
    }
  }{ L^q ( P_{\gm_s} ) } 
  - 
  L^* 
  \abra{
    \Norm{ 
      \abs{ 
        \nab_{d} Q_s f 
      }
    }{ L^q ( P_{\gm_s} ) } 
  } 
} 
d s 
\\
& \le 
L ( \tilde{d} (x,y) )
,
\end{align*}
where the second inequality 
comes from the definition of $L^*$ 
as the Legendre conjugate. 
Substituting this estimate into \eq{Kantorovich2}, 
we obtain the desired estimate.  

\textbf{(iii) The case $p=\infty$:} 
Since \eq{gradient} holds with $q=1$, 
the H\"{o}lder inequality implies 
\eq{gradient} for any $q > 1$. 
Thus we obtain \eq{cost} for any $1 \le p < \infty$. 
Therefore, 
by virtue of \Lem{infinite-cost}, 
the conclusion follows 
by tending $p$ to $\infty$ in \eq{cost}. 
\epf
\brem{Orlicz}
Our duality between $L^p$ and $L^q$ 
can be extended to 
a similar one between Orlicz norms. 
In fact, 
there are H\"{o}lder-type inequalities 
(see \cite{Adams-Fourn}, for instance) 
which will be used in the implication 
\ref{transport-bound} 
$\Rightarrow$ 
\ref{gradient-estimate}. 
For the converse, 
all properties of Hamilton-Jacobi semigroup 
we will use in the proof still hold 
in such a generality. 
\erem 
\brem{support}
If \eq{cost} holds with $p > 1$, 
then we obtain 
the following slightly stronger version of 
($G_\infty$); 
for any $f \in C_{b,L} (X)$ and $x \in X$, 
\begin{equation} 
\tag{$G_\infty'$} 
\label{eq:gradient-infinity-minus}
\abs{ \nab_{\tilde{d}} P f } (x) 
\le 
\Norm{ \abs{ \nab_{d} f } }{ L^\infty ( P_x ) } 
. 
\end{equation}
As we have seen in the proof of \Prop{gradient-cost}, 
a weaker condition ($G_\infty$) is sufficient 
to obtain ($C_1$). 
At this moment, 
the author does not know any example 
that ($C_p$) holds only for $p=1$ 
and ($G_\infty'$) fails. 
\erem
\section{Applications}
\label{sec:application} 

In a class of sub-Riemannian manifolds, 
$L^q$-gradient estimates of a subelliptic heat semigroup 
is shown recently by an analytic method. 
In these cases, 
we can obtain the corresponding 
$L^p$-Wasserstein control 
via \Thm{duality} 
though their notion of gradient 
looks different from ours. 
To explain how we deal with it, 
we will demonstrate 
a general framework of sub-Riemannian geometry 
generated by a family of vector fields. 
We refer to 
\cite{SmP,Mont-sRiem,Strich_subRiem} 
for details. 

Throughout this section, we assume $X$ to be 
a finite dimensional, $\sigma$-compact, connected, 
smooth differentiable manifold. 
Consider a family of vector fields 
$\{ X_1 , \cdots , X_n \}$ on $X$. 
We assume that 
$\{ X_i (x) \}_{i=1}^n$ is linearly independent 
on $T_x X$ for all $x \in X$ and 
that $\{ X_i \}_{i=1}^n$ satisfies 
the H\"{o}rmander condition. 
The latter one means that there exists a number $m$ 
such that the family of vector fields 
generated by $\{ X_i \}_{i=1}^n$ and 
their commutators up to the length $m$ 
spans $T_x X$ for each $x \in X$. 
Let $\cH \subset T X$ be the subbundle 
generated by $\{ X_i \}_{i=1}^n$; 
$
\cH_x 
: = 
\Span 
\bbra{ 
  X_1 (x) , \ldots , X_n (x) 
}
$. 
We define a metric on $\cH$ such that 
$\{ X_i (x) \}_{i=1}^n$ becomes 
an orthonormal basis of $\cH_x$ for $x \in X$. 
We are interested in the case $\cH \neq TX$. 
Associated with this metric, 
we define a function $d$ on $X$ as follows. 
We say a piecewise smooth curve 
$\gm \: : \: [0, l] \to X$ horizontal 
if $\dot{\gm}(t) \in \cH_{\gm (t)}$ for every $t$ 
where $\gm$ is differentiable. 
For $x , y \in X$, we define $d(x,y)$ by 
\[
d(x,y) 
: = 
\inf \bbra{ 
  \left. 
  \int_0^l 
  \Norm{\dot{\gm} (t) }{\cH_{\gm (t)} } 
  dt 
  \;
  \right| 
  \begin{array}{l}
  \gm \: : \: [0,l] \to X 
  \mbox{ horizontal curve}, 
  \\ 
  \gm (0) = x 
  , \, 
  \gm (l) = y 
  \end{array}
}. 
\] 
By the Chow theorem, 
the H\"{o}rmander condition ensures 
that $d(x,y) < \infty$ for $x,y \in X$. 
As a result, 
the function 
$d \: : \: X \times X \to [0,\infty)$ 
becomes a distance. 
It is called the Carnot-Caratheodory distance. 
Note that the topology determined by $d$ coincides 
with the original one on $X$. 
We assume that $( X , d )$ is complete. 

Let $v$ be a Borel measure on $X$ such that 
its restriction on each local coordinate 
has a smooth density with respect to 
the Lebesgue measure associated with the coordinate. 
Let $\Dl_{\cH} : = \sum_{i=1}^n X_i^* X_i / 2$ be 
the sub-Laplacian 
associated with $\{ X_i \}_{i=1}^n$ and $v$. 
Here $X_i^*$ is the adjoint operator of $X_i$ 
with respect to $v$. 
By the completeness of $d$, 
$\Dl_\cH$ is essentially selfadjoint 
(see \cite{Strich_subRiem}). 
Take the selfadjoint extension of $\Dl_\cH$ 
(also denoted by $\Dl_\cH$) 
and consider the associated heat semigroup 
$P_t = \exp ( t \Dl_{\cH} / 2 )$. 
By the hypoellipticity of $\Dl_\cH$, 
$P_t$ has a smooth density function 
with respect to $v$. 
In particular, 
$P_t$ becomes a Feller semigroup. 
We assume that $P_t$ is conservative, i.e. $P_t 1 = 1$. 
For a smooth function $f \: : \: X \to \R$, 
we define the carr\'{e} du champ operator 
$\Gm (f) \: : \: X \to \R$ 
by 
$
\Gm (f) (x) 
= 
\sum_{i=1}^n \abs{ X_i f (x) }^2 
$. 

An $L^q$-gradient estimate for $P_t$ 
associated with $\Gm$ 
is formulated as follows; 
given $q \in [1,\infty)$, 
there exists $K_q (t) > 0$ for each $t > 0$ 
such that, 
for any $f \in C^\infty_c (X)$, 
\begin{equation} 
\label{eq:smooth-gradient}
\Gm ( P_t f )(x)^{1/2} 
\le 
K_q (t) 
\bbra{ 
  P_t 
  \abra{ 
    \Gm ( f )^{q/2} 
  } (x) 
}^{1/q} , 
\end{equation}
where $C^\infty_c (X)$ is the set of all smooth functions 
$f \: : \: X \to \R$ with compact supports. 
As we see in the following, 
\eq{smooth-gradient} 
implies our gradient estimate. 
\bprop{approximation} 
\eq{smooth-gradient} 
for $f \in C^\infty_c (X)$ 
implies 
\eq{gradient} 
for $P = P_t$, $\tilde{d} = K_q (t) d$ and 
any $f \in C_{L} (X)$ with a compact support. 
\eprop 
\bpf
First we extend \eq{smooth-gradient} 
for $f \in C_{b,L}(X)$. 
By virtue of Corollary~11.8 of \cite{SmP}, 
for $f \in C_{b,L} (X)$, 
the distributional derivatives 
$\{ X_i f \}_{i=1}^n$ are represented 
as a bounded functions and 
$\abs{ \Gm f }^{1/2} \le \Norm{\nab_d f}{\infty}$ 
holds $v$-almost everywhere. 
Moreover, Theorem~11.7 of \cite{SmP} implies 
$\abs{ \Gm f }^{1/2} \le g_f$ for any upper gradient $g_f$. 
In particular, \Lem{upper-gradient} implies 
$\abs{ \Gm f }^{1/2} \le \abs{ \nab_d f }$. 
Though they discussed 
the case that 
$X$ is an open subset of a Euclidean space 
in \cite{SmP}, 
we can extend it to our case 
with the aid of a partition of unity. 
By a mollifier argument 
together with use of a partition of unity again, 
we can take a sequence 
$f_k \in C^\infty_c (X)$ 
such that $f_k \to f$ and $\Gm f_k \to \Gm f$ 
almost surely (cf. \cite[Theorem~11.9]{SmP}). 
Thus \eq{smooth-gradient} holds 
for any $f \in C_{b,L} (X)$ with a compact support. 

Note that 
$\abs{\Gm f }^{1/2}$ is an upper gradient 
if $f \in C^\infty (X)$ 
(see \cite[Proposition~11.6]{SmP}, for instance). 
Since $P_t f \in C^\infty (X)$ in our case, 
for a minimal geodesic $\gm$ joining $x$ and $y$, 
\begin{align*}
P_t f (x) - P_t f (y) 
& \le 
\int_0^{d(x,y)} 
\bbra{ 
  \Gm ( P_t f )( \gm (s) ) 
}^{1/2} ds 
\\
& \le 
\int_0^{d(x,y)} 
\bbra{ 
  P_t \abra{
    \Gm ( f )^{q/2} 
  }
  ( \gm (s) ) 
}^{1/q} ds 
\\
& \le 
\int_0^{d(x,y)} 
\bbra{ 
  P_t \abra{
    \abs{ \nab_d f }^q 
  }
  ( \gm (s) ) 
}^{1/q} ds 
. 
\end{align*}
Hence the conclusion follows 
by dividing the above inequality by $d(x,y)$ 
and by letting $y \to x$.  
\epf 
\brem{upper-gradient}
If we suppose 
\Ass{regularity} \ref{doubling} \ref{Poincare} 
in \Prop{approximation}, 
then Theorem 6.1 of \cite{Chee99} asserts 
that the minimal generalized upper gradient of $f$ 
coincides with $\abs{ \nab f }$ almost everywhere. 
Since the first part of 
the proof of \Prop{approximation} 
implies that $\abs{ \Gm f }^{1/2}$ is 
the minimal generalized upper gradient 
for $f \in C_{L} (X)$ with a compact support, 
the proof can be completed there in this case. 
\erem 

As far as the author knows, 
\eq{smooth-gradient} is established 
in the following cases;  
\begin{itemize}
\item 
The case $q=1$ with $K_1 (t) \equiv K$ for some $K > 0$ 
on groups of type H \cite{Eldr09} 
(including the Heisenberg group of arbitrary dimension, 
see \cite{BBBC08,HQLi06} also). 
\item
The case $q > 1$ on an arbitrary Lie group 
\cite{Mel08}. 
Especially, $K_p (t) \equiv K_p$ for some $K_p > 0$ 
if it is nilpotent. 
\item
The case $q > 1$ with $K_q (t) = K_q \e^{-t}$ 
for some $K_q > 0$ 
on $\mathbf{SU}(2)$ \cite{BauBonn08}. 
\end{itemize}
In all these cases, $v$ is chosen 
to be a right-invariant Haar measure 
and hence the associated sub-Laplacian 
is of the form $\Dl_{\cH} = \sum_{i=1}^n X_i^2$. 
All conditions in \Ass{regularity} 
hold in these cases. 
For \ref{density}, we have already observed. 
By the homogeneity of the space, 
we can reduce the assertion 
in the case of a Euclidean domain 
(see \Rem{Poincare} also). 
Thus 
\ref{doubling} 
and 
\ref{Poincare} with $p_0 = 1$ 
follow 
from 
Theorem 11.19 
and  
Theorem 11.21 of \cite{SmP}. 
Note that 
\eq{smooth-gradient} is shown 
on a wider class of functions 
than $C^\infty_c (X)$ in some cases. 
But it is not necessary for our purpose. 

Combining \Prop{approximation} with 
\Thm{duality} in these cases, 
we obtain \eq{cost} 
for $P = P_t$ and $\tilde{d} = K_q (t) d$. 
Though 
$f$ is restricted to have a compact support 
in \Prop{approximation}, 
it is sufficient to show \eq{cost} 
(see \Rem{Lipschitz-approximation}). 

The following simple examples explain 
a probabilistic meaning of these consequences. 
\bexam{Heisenberg}
The 3-dimensional Heisenberg group is realized on $\R^3$ 
with the multiplication defined by 
\[
( x , y , z )
\cdot 
( x' , y' , z' ) 
= 
\abra{ 
  x + x' , 
  y + y' , 
  z + z' + \frac12 ( x y' - y x' )
}. 
\]
The Lebesgue measure $v$ on $\R^3$ is 
a bi-invariant Haar measure. 
Let us define left-invariant vector fields 
$X, Y$ and $Z$ by 
\[
X := \pdel{}{x} - \frac{y}{2} \pdel{}{z}, 
\qquad 
Y := \pdel{}{y} + \frac{x}{2} \pdel{}{z}, 
\qquad  
Z := \pdel{}{z}. 
\]
Set $\cH := \Span \{ X , Y \}$. 
Then the diffusion process 
$\{ \mathbf{B}_t^{\mathbf{x}} \}_{t \ge 0}$ 
associated with 
$\Dl_{\cH} / 2 = ( X^2 + Y^2 )/2$ 
starting at 
$\mathbf{x} = (x,y,z) \in \R^3$ 
is given by 
\[
\mathbf{B}_t^{\mathbf{x}} 
: = 
\abra{ 
  x + W^{(1)}_t , 
  y + W^{(2)}_t , 
  z 
  + \frac12 
  \int_0^t 
  ( x + W^{(1)}_s ) d W^{(2)}_s 
  - ( y + W^{(2)}_s ) d W^{(1)}_s 
}
, 
\]
where $( W^{(1)}_t , W^{(2)}_t )$ is 
a Brownian motion on $\R^2$. 
It means that the diffusion process associated 
with $\Dl_{\cH}/2$ is given 
by the 2-dimensional Euclidean Brownian motion and 
the associated L\'{e}vy stochastic area. 
The corresponding heat semigroup is given by 
$
P_t f ( \mathbf{x} ) 
= 
\E \cbra{ 
  f ( \mathbf{B}_t^{\mathbf{x}} ) 
}
$ for $f  \in C_b (X)$. 
In this framework, \eq{smooth-gradient} 
for $q=1$, $P = P_t$ and $K_1 (t) \equiv K$ 
is shown in \cite{BBBC08,HQLi06}. 
Thus we obtain ($C_\infty$). 
It means that,   
for each $t > 0$ and 
$\mathbf{x} , \mathbf{y} \in \R^3$, 
there exists a coupling 
$
( 
\bar{\mathbf{B}}^{\mathbf{x}}_t 
, 
\bar{\mathbf{B}}^{\mathbf{y}}_t 
)
$ 
of 
$\mathbf{B}_t^{\mathbf{x}}$ 
and 
$\mathbf{B}_t^{\mathbf{y}}$ 
such that  
\begin{equation} 
\label{eq:Heisenberg-bound}
d \abra{ 
  \bar{\mathbf{B}}_t^{\mathbf{x}} , 
  \bar{\mathbf{B}}_t^{\mathbf{y}} 
} 
\le 
K d ( \mathbf{x} , \mathbf{y} )
\end{equation} 
holds almost surely. 
Here $d$ is the Carnot-Caratheodory distance 
associated with $\cH$. 
In this case, it is known that 
$d$ is equivalent to the so-called Kor\'{a}nyi distance. 
That is, there exist constants $C_1 , C_2 > 0$ 
such that, for any 
$
\mathbf{x} = ( x , y , z ) , 
\mathbf{y} = ( x' , y' , z' ) 
\in \R^3
$, 
\[
C_1 d ( \mathbf{x},\mathbf{y} )
\le 
\bbra{ 
  \abra{ 
    ( x - x' )^2 
    + 
    ( y - y' )^2 
  }^2 
  + 
  \abra{ 
    z - z' + \frac12 ( x y' - y x' ) 
  }^2 
}^{1/4}
\le 
C_2 d ( \mathbf{x} , \mathbf{y} ) . 
\] 
Thus \eq{Heisenberg-bound} is 
also interpreted 
in terms of the Kor\'{a}nyi distance. 
\eexam
\brem{pathwise}
In \Exam{Heisenberg}, 
($C_\infty$) provides only a coupling of 
$\mathbf{B}^{\mathbf{x}}_t$ and $\mathbf{B}^{\mathbf{y}}_t$ 
for each fixed $t>0$. 
%
When $X$ is a Riemannian manifold, 
($C_\infty$) holds 
if and only if 
there exists a coupling 
$( 
\bar{\mathbf{B}}^{\mathbf{x}}_t 
, 
\bar{\mathbf{B}}^{\mathbf{y}}_t 
)_{t \ge 0}$
of two Brownian motions 
$( \mathbf{B}^{\mathbf{x}}_t )_{t \ge 0}$ and 
$( \mathbf{B}^{\mathbf{y}}_t )_{t \ge 0}$ 
starting from $\mathbf{x}$ and $\mathbf{y}$ 
respectively such that \eq{Heisenberg-bound} 
holds 
for every $t \ge 0$ with $K = \e^{-kt}$ almost surely 
(see \cite{Stu_Renes05}, for instance). 
In \Exam{Heisenberg}, 
it is not clear whether 
a similar result holds or not. 
Actually, in Riemannian case, 
the fact that 
the constant $\e^{-kt}$ is multiplicative 
in $t \ge 0$ 
plays a prominent role 
to construct a coupling of Brownian motions 
from a control of their infinitesimal motions. 
As observed in \cite{Dri-Mel05}, 
we cannot expect such a multiplicativity 
in the case of \Exam{Heisenberg}. 
\erem
\bexam{Levy-areas}
On $\R^n \times \R^{n(n-1)/2}$, 
we introduce 
a structure of nilpotent Lie group of step 2 
as follows; 
for 
$
\mathbf{x} 
= 
( 
  ( x_i )_{i=1}^n 
  ; 
  ( z_{ij} )_{1 \le i < j \le n} 
), 
\mathbf{y} 
= 
( 
  ( x_i' )_{i=1}^n 
  ; 
  ( z_{ij}' )_{1 \le i < j \le n} 
) 
\in \R^n \times \R^{n(n-1)/2} 
$, 
\[
\mathbf{x} \cdot \mathbf{y} 
= 
\abra{ 
  ( x_i + x_i' )_{i=1}^n ; 
  \abra{ 
    z_{ij} + z_{ij}' 
    + \frac12 
    ( x_i x_j' - x_j x_i' ) 
  }_{1 \le i < j \le n}
}
. 
\]
As in \Exam{Heisenberg}, 
the Lebesgue measure $v$ 
on $\R^n \times \R^{n(n-1)/2}$ 
becomes a bi-invariant Haar measure. 
Let us define left-invariant vector fields 
$\{ X_i \}_{i=1}^n$ and $\{ Z_{ij} \}_{1 \le i < j \le n}$ 
by 
\begin{align*} 
X_i 
:= 
\pdel{}{x_i} 
- 
\sum_{i < j \le n} \frac{x_j}{2} \pdel{}{z_{ji}} 
+ 
\sum_{1 \le j < i} \frac{x_j}{2} \pdel{}{z_{ij}} , 
\qquad
Z_{ij} 
:= 
\pdel{}{z_{ij}} 
.
\end{align*} 
Set $\cH : = \Span \{ X_i \}_{i=1}^n$. 
The diffusion process 
$\{ \mathbf{B}_t^{\mathbf{x}} \}_{t \ge 0}$ 
associated with the sub-Laplacian 
$\Dl_{\cH}/2 = \sum_{i=1}^n X_i^2 /2$ 
starting at 
$
\mathbf{x} 
= 
( 
\{ x_i \}_{i=1}^n ; 
\{ z_{ij} \}_{i=1}^n 
) 
\in \R^n \times \R^{n(n-1)/2}
$
is given by 
\[
\mathbf{B}^{\mathbf{x}}_t 
 = 
\abra{ 
  \abra{ 
    x_i + W^{(i)}_t 
  }_{i=1}^n ; 
  \abra{ 
    z_{ij} 
    + \frac12 
    \int_0^t 
    ( x_i + W^{(i)}_s ) d W^{(j)}_s 
    - 
    ( x_j + W^{(j)}_s ) d W^{(i)}_s  
  }_{1 \le i < j \le n }
}
.
\]
We can easily verify that 
this group is of type H 
only if $n = 1$ 
(see Corollary 1 of \cite{Kaplan_typeH}, for example).  
But it is still in the framework of \cite{Mel08}. 
Thus, for each $p \in [ 1 , \infty )$, 
there is a constant $K_p > 0$ 
such that, 
for any pair 
$\mathbf{x}, \mathbf{y} \in \R^n \times \R^{n(n-1)/2}$, 
there is a coupling 
$
( 
\bar{\mathbf{B}}^{\mathbf{x}}_t 
, 
\bar{\mathbf{B}}^{\mathbf{y}}_t 
)
$ 
of 
$\mathbf{B}_t^{\mathbf{x}}$ 
and 
$\mathbf{B}_t^{\mathbf{y}}$ 
satisfying 
\begin{equation} 
\label{eq:Heisenberg-bound2}
\E \cbra{ 
  d 
  \abra{ 
    \bar{\mathbf{B}}_t^{\mathbf{x}} , 
    \bar{\mathbf{B}}_t^{\mathbf{y}} 
  }^p  
}^{1/p}
\le 
K_p d ( \mathbf{x} , \mathbf{y} ) . 
\end{equation} 
Finally, we give a remark that 
a different kind of coupling of 
this process is studied 
by Kendall \cite{Kend_Levy-couple}. 
He showed the existence of 
a successful coupling. 
As mentioned there, 
studying a coupling of this process 
has a possibility of 
a future application 
to rough path theory \cite{Friz-Victoir,Lyons-Qian}. 
\eexam 

\providecommand{\bysame}{\leavevmode\hbox to3em{\hrulefill}\thinspace}
\providecommand{\MR}{\relax\ifhmode\unskip\space\fi MR }
\providecommand{\MRhref}[2]{%
  \href{http://www.ams.org/mathscinet-getitem?mr=#1}{#2}
}
\providecommand{\href}[2]{#2}

\vspace{1cm}
\begin{flushright}
Kazumasa Kuwada\\
\vspace{.3cm}
\small
Graduate School of Humanities and Sciences 
\\
\small
Ochanomizu University \\
\small
Tokyo 112-8610, Japan
\\
\vspace{.2cm}
\small
\textit{e-mail}: \texttt{kuwada@math.ocha.ac.jp}
\end{flushright}
\end{document}